\documentclass [6pt,a4paper]{article}
\usepackage{bbm}
\usepackage{amsfonts}
\usepackage{mathrsfs}
\usepackage {amssymb}
\usepackage {amsmath}
\usepackage{amsthm}
\usepackage{latexsym}
\def\dse#1{\vskip 0.6cm\noindent
        {\large\bf #1}
        \vskip 0.4cm}
\newcommand{\pf}{{\bf Proof. \ }}
\def\dse#1{\vskip 0.6cm\noindent
        {\large\bf #1}
        \vskip 0.4cm}
\usepackage{lastpage}
\usepackage{epsfig}

\begin{document}
\begin{center}
\textbf{\large{ One-Lee weight and two-Lee weight
$\mathbb{Z}_2\mathbb{Z}_2[u]$-additive codes}}\footnote{E-mail
addresses: luzhenliang1992@sina.cn(Z.Lu),
liqiwangg@163.com(L.Wang), zhushixin@hfut.edu.cn(S.Zhu),
kxs6@sina.com(X.Kai).}
\end{center}

\begin{center}
{ { Zhenliang Lu, \ Liqi Wang, \ Shixin Zhu, \ Xiaoshan Kai } }
\end{center}

\begin{center}
\textit{\footnotesize School of Mathematics, Hefei University of
Technology, Hefei 230009, Anhui, PR China}
\end{center}

\noindent\textbf{Abstract} In this paper, we study one-Lee weight
and two-Lee weight codes over $\mathbb{Z}_{2}\mathbb{Z}_{2}[u]$,
where $u^{2}=0$. Some properties of one-Lee weight
$\mathbb{Z}_{2}\mathbb{Z}_{2}[u]$-additive codes are given, and a
complete classification of one-Lee weight
$\mathbb{Z}_2\mathbb{Z}_2[u]$-additive formally self-dual codes is
obtained. The structure of two-Lee weight projective
$\mathbb{Z}_2\mathbb{Z}_2[u]$ codes is determined. Some optimal
binary linear codes are obtained directly from one-Lee weight and
two-Lee weight $\mathbb{Z}_{2}\mathbb{Z}_{2}[u]$-additive codes via
the extended Gray map.\\


\dse{1~~Introduction} A code over a finite field is said to be a
$t$-weight code if the cardinality of the set of nonzero Hamming
weights is $t$. This class of codes is interesting in the area of
coding theory due to many applications in mathematics such as
strongly regular graphs and projective point-sets ([2,4,5,7,9]).
Thus, $t$-weight codes over finite fields have received a lot of
attention (see [8,13,15,16]). After a landmark paper [14] showed
that certain good nonlinear binary codes can be constructed from
linear codes over $\mathbb{Z}_{4}$ via the Gray map, $t$-Lee weight
codes over finite rings have been extensively studied (see, for
example, [1,6,11,17,18,20]). One-Lee weight linear codes over
$\mathbb{Z}_{4}$ were first studied by Carlet [6]. More generally,
the structure of one-Lee weight linear codes over $\mathbb{Z}_{m}$
was determined by Wood [21]. The structures of one-Lee weight codes
and two-Lee weight projective codes over
$\mathbb{F}_2+u\mathbb{F}_2$ were explored by Wang et al. [20]. The
results [20] were extended to the ring $\mathbb{F}_p+u\mathbb{F}_p$
by Shi et al. and many optimal linear codes over $\mathbb{F}_{p}$
were obtained from one-Lee weight codes and two-Lee weight codes
over $\mathbb{F}_p+u\mathbb{F}_p$ [17,18].

Additive codes were first introduced by Delsarte [10], where it was
shown that any abelian binary propelinear code has the form
$\mathbb{Z}_{2}^{\alpha}\times\mathbb{Z}_{4}^{\beta}$. This
motivated one to study $\mathbb{Z}_{2}\mathbb{Z}_{4}$-additive
codes. The structure of one weight
$\mathbb{Z}_{2}\mathbb{Z}_{4}$-additive codes has been investigated
by Dougherty et al. [12]. Formally self-dual codes in
$\mathbb{Z}_{2}^{\alpha}\times\mathbb{Z}_{4}^{\beta}$ are
constructed and some good formally self-dual codes as the binary
images of these codes were obtained [11,22]. Recently,
$\mathbb{Z}_{2}\mathbb{Z}_{2}[u]$-additive codes have been
introduced by Aydogdu et al. [1], where the standard forms of
generator and parity check matrices of the $\mathbb Z_{2}\mathbb
Z_{2}[u]$-additive codes were determined and a MacWilliams identity
for these codes was proved.

Motivated by the work listed above, the purpose of this paper is to
discuss one-Lee weight and two-Lee weight codes over
$\mathbb{Z}_{2}\mathbb{Z}_{2}[u]$. We study the structure and
possible weights for all one-Lee weight
$\mathbb{Z}_2\mathbb{Z}_2[u]$-additive formally self-dual codes. In
particular,
the structure of two-Lee weight projective $\mathbb{Z}_2\mathbb{Z}_2[u]$-additive codes are discussed.
Some optimal linear codes as the binary images of one-Lee and two-Lee weight $\mathbb{Z}_2\mathbb{Z}_2[u]$-additive codes are given.\\

\dse{2~~Preliminaries}

Denote by $\mathbb{Z}_2$ the ring of integers module 2. Let $R$ be
the commutative ring $\mathbb{Z}_2+u\mathbb{Z}_2=\{a+ub\mid a,b\in
\mathbb{Z}_2\}=\{0, 1, u, u+1\}$, where $u^2=0$. Obviously, the ring
$\mathbb{Z}_2$ is a subring of the ring $R$. Let $\mathbb{Z}_2^{n}$
be the set of $n$-tuples over $\mathbb{Z}_2$. Any non-empty subset
$C$ of $\mathbb{Z}_2^{n}$ is called a binary code, and if that code
is a vector space, then the code $C$ is called a
$\mathbb{Z}_2$-linear code. We define a Gray map $\psi$ from $R$ to
$\mathbb{Z}_{2}^{2}$ in the following way:
\begin{align*}
\psi:R&\rightarrow \mathbb{Z}_{2}^{2},\\
 (a+ub)&\rightarrow(b, a+b).
\end{align*}
So $\psi(0)=(0,0), \psi(1)=(0,1), \psi(u)=(1,1), \psi(u+1)=(1,0).$
The ring $\mathbb{Z}_2\mathbb{Z}_2[u]$ is defined by
$$\mathbb{Z}_2\mathbb{Z}_2[u]=\{(a|b)|a\in\mathbb{Z}_2 \ \textrm{and} \ b\in R\}.$$

The ring is not closed under standard multiplication by the element
$u$ in the ring $R$. Therefore, to make the ring
$\mathbb{Z}_2\mathbb{Z}_2[u]$ an $R$-module, we introduce the
following scalar multiplication:

 for any $c_{1}=({a_{0},a_{1},\cdots,a_{\alpha-1}|b_{0},b_{1},\cdots,b_{\beta-1}})\in\mathbb{Z}_{2}^{\alpha}R^{\beta}$
and $c=r+qu\in R,$
$$c\cdot c_{1}=({ra_{0},ra_{1},\cdots,ra_{\alpha-1}|cb_{0},cb_{1},\cdots,cb_{\beta-1}}).$$

Let $(a|b)\in\mathbb{Z}_2^{\alpha}R^{\beta}$,
where $\alpha$, $\beta$ are nonnegative integers, $
a=({a_{0},a_{1},\cdots,a_{\alpha-1}})\in \mathbb{Z}_2^{\alpha} $,
and $ b=(b_{0},b_{1},\cdots,b_{\beta-1})\in R^{\beta}$. Then the
Gray map $\Phi$ from
$\mathbb{Z}_2^{\alpha}R^{\beta}$ to
$\mathbb{Z}_2^{\alpha+2\beta}$ is defined as
$$\Phi(a|b)=({a_{0},a_{1},\cdots,a_{\alpha-1}|\psi(b_{0}),\psi(b_{1}),\cdots,\psi(b_{\beta-1})}).$$

A code $C$ is called a $\mathbb{Z}_2\mathbb{Z}_2[u]$ additive code
if it is a $R$-submodule of
$\mathbb{Z}_2^{\alpha}R^{\beta}$ under the above
scalar multiplication. For any two vectors
$v=(v_{1},v_{2},\ldots,v_{\alpha}\mid
v_{\alpha+1},\ldots,v_{\alpha+\beta})$ and
$w=(w_{1},w_{2},\ldots,w_{\alpha}\mid
w_{\alpha+1},\ldots,w_{\alpha+\beta})\in
$$\mathbb{Z}_2^{\alpha}R^{\beta}$, their inner product is
defined by
\[v\cdot w=u \sum_{i=1}^{\alpha}v_{i}w_{i}+\sum_{j=\alpha+1}^{\alpha+\beta}v_{j}w_{j}\in R.\]
The dual code of $C$ is defined as
\[C^{\perp}=\{w\in\mathbb{Z}_2^{\alpha}R^{\beta}\mid
v\cdot w=0~ \textrm{for}~\textrm{all}~v\in C\}.\]

For any codeword $ c=(a_{0},a_{1},\cdots,a_{\alpha-1}|$
$b_{0},b_{1},\cdots,b_{\beta-1})\in C $, we denote by $w_{L}(c)$ the
Lee weight in the following way:
$w_{L}(c)=\sum_{i=0}^{\alpha-1}w_{L}(a_{i})+\sum_{i=0}^{\beta-1}w_{L}(b_{i}),$
where $w_{L}(0)=0, w_{L}(1)=1, w_{L}(u)=2$ and $w_{L}(u+1)=1$. Note
the Hamming weight of an $n$-tuple is the number of its nonzero
entries. So, $w_{L}(a_{i})=w_{H}(a_{i})$ for $0\leq i\leq\alpha-1$,
and $w_{L}(b_{i})=w_{H}(\psi(b_{i}))$ for $0\leq i\leq\beta-1$.
Therefore,
$w_{L}(c)=\sum_{i=0}^{\alpha-1}w_{H}(a_{i})+\sum_{i=0}^{\beta-1}w_{H}(\psi(b_{i})).$
This Lee weight function also defines a Lee distance function
$d_{L}(x,y)=w_{L}(x-y)$ on $\mathbb{Z}_2\mathbb{Z}_2[u]$ between
$x,y\in\mathbb{Z}_2^{\alpha}R^{\beta}$. It is
easy to see that the Gray map $\Phi$ is a weight preserving map from
$(\mathbb{Z}_2^{\alpha}R^{\beta},~\textrm{Lee}~~\textrm{weight})
~\textrm{to}~(\mathbb{Z}_2^{\alpha+2\beta},$
$~\textrm{Hamming}~~\textrm{weight})$. And $\Phi$ is a distance
preserving map from
$(\mathbb{Z}_2^{\alpha}R^{\beta},~\textrm{Lee}~~
\textrm{distance})~\textrm{to}~(\mathbb{Z}_2^{\alpha+2\beta},~\textrm{Hamming}~~\textrm{distance})$.\\

\noindent\textbf{Theorem 2.1} ([1]). \emph{Let $C$ be a
$\mathbb{Z}_2\mathbb{Z}_2[u]$-additive code. Then $C$ is permutation
equivalent to a $\mathbb{Z}_2\mathbb{Z}_2[u]$-additive code with the
standard form matrix
$$ \begin {pmatrix} I_{k_{0}} & A_{1} & 0 & 0 & uT\\0 & S & I_{k_1} & A &  B_{1}+uB_{2} \\ 0 & 0 & 0 & uI_{k_{2}} & uD\end{pmatrix},$$
where $A, A_{1}, S, D, B_{1}, B_{2}$ and T are
$\mathbb{Z}_2$-matrices, $I_{k_i}(i=0,1,2)$ are $k_{i}\times k_{i}$  identity matrices.
The code $C$ contains $2^{k_{0}+2k_{1}+k_{2}}$
codewords.} \\

Related to the standard form matrix, such an additive
code $C$ is isomorphic to
$\mathbb{Z}_2^{k_{0}}\times\mathbb{Z}_2^{2k_{1}}\times\mathbb{Z}_2^{k_{2}}$,
and $C$ is said to be of type $(\alpha, \beta; k_{0}, k_{1},
k_{2})$. Throughout this paper, we make the convention that the
standard form matrix of $C$ has no all zero columns.

\dse{3~~One-Lee weight $\mathbb{Z}_2\mathbb{Z}_2[u]$-additive codes}

In this section, we study the properties of one-Lee weight $\mathbb{Z}_2\mathbb{Z}_2[u]$-additive codes. We first give the following definition.\\

\noindent\textbf{Definition 3.1.} A linear code is said to have one
weight if every nonzero codeword has the same Hamming weight. A
$\mathbb{Z}_2\mathbb{Z}_2[u]$-additive code is called a one-Lee
weight code if every nonzero codeword has the same Lee weight.\\

 Let
$C$ be a one-Lee weight $\mathbb{Z}_2\mathbb{Z}_2[u]$-additive codes
with nonzero weight $m$. It is known that $\Phi(C)$ is a one weight
binary code
with weight $m$.\\

\noindent\textbf{Theorem 3.2.}~\emph{Let $C$ be a one-Lee weight
$\mathbb{Z}_2^{\alpha}R^{\beta}$-additive code, then
the sum $ \sum_{c\in C}w_{L}(c)$ of the Lee weights of all the
codewords of $C$ is equal to $\frac{|C|}{2}(\alpha+2\beta)$.}\\

\noindent\pf Write all the codewords of $C$ as a
$|C|\times(\alpha+\beta)$ array. Note that $C$ is an additive code.
First, we consider the binary part of the codewords, it is easy to
check that the number of each column contains 0 and 1 is equally.
Next, we consider the nonbinary part, namely the last $\beta$
columns of array. Clearly, any column in this part corresponds to
the next two cases and does not contain any other case.
\begin{itemize}
\item[(1)] Any column contains an equal numbers of $0$, $1$, $u$ and $u+1$.
\item[(2)] Any column contains an equal numbers of $0$ and $u$.
\end{itemize}
Let $N_{1}$ be the number of columns in Case (1). Then $\beta-N_{1}$
is the number of columns in Case (2). Hence, we have

\begin{align*}
\sum_{c\in C}w_{L}(c)=&\alpha\left(\frac{|C|}{2}\right)+N_{1}\left(\frac{|C|}{4}+\frac{|C|}{2}\cdot2+\frac{|C|}{4}\right)+(\beta-N_{1})\left(\frac{|C|}{2}\cdot2\right)\\
=&\frac{|C|}{2}(\alpha+2\beta).
\end{align*}\qed\\

\noindent\textbf{Theorem 3.3.}~\emph{ Let $C$ be a one-Lee weight $\mathbb{Z}_2\mathbb{Z}_2[u]$-additive
code of type $(\alpha, \beta; k_{0}, k_{1}, k_{2})$ with weight $m$, then there exists a unique
positive integer $\lambda$ such that $m=\lambda\frac{|C|}{2}$ and $\alpha+2\beta=\lambda(|C|-1)$.}\\

\noindent\pf According to Theorem 3.2, the sum $ \sum_{c\in
C}w_{L}(c)$ of the Lee weights of all the codewords of $C$ is equal
to $\frac{|C|}{2}(\alpha+2\beta)$. On the other hand, the sum of the
weights
is equal to $m(|C|-1)$. So, $\frac{|C|}{2}(\alpha+2\beta)=m(|C|-1)$. Since $|C|=2^{k_{0}+2k_{1}+k_{2}}$, it follows that $\gcd(\frac{|C|}{2},|C|-1)=1$. Hence, there must exists a unique positive integer $\lambda$ such that $m=\lambda\frac{|C|}{2}$ and $\alpha+2\beta=\lambda(|C|-1)$.\qed\\

\noindent\textbf{Corollary 3.4.}~\emph{Let $C$ be a one-Lee weight $\mathbb{Z}_2\mathbb{Z}_2[u]$-additive code
of type $(\alpha, \beta; k_{0}, k_{1}, k_{2})$ with weight $m$. If $m$ is odd,
then $C=\{(\underbrace{0,\cdots,0}_{\alpha}|\underbrace{0,\cdots,0}_{\beta}), (\underbrace{1,\cdots,1}_{\alpha}|\underbrace{u,\cdots,u}_{\beta})\}$.} \\\\
\noindent\pf By Theorem 3.3, $m=\lambda\frac{|C|}{2}$ and
$\alpha+2\beta=\lambda(|C|-1)$.
By Theorem 2.1, $|C|=2^{k_{0}+2k_{1}+k_{2}}$. Since $m$ is odd, we can deduce that $\lambda$ is odd and $\frac{|C|}{2}=1$.
Then $\alpha+2\beta=\lambda=m$. Notice that $(\underbrace{1,\cdots,1}_{\alpha}|\underbrace{u,\cdots,u}_{\beta})$ is
the only codeword with weight $\alpha+2\beta$, so
$$C=\{(\underbrace{0,\cdots,0}_{\alpha}|\underbrace{0,\cdots,0}_{\beta}), (\underbrace{1,\cdots,1}_{\alpha}|\underbrace{u,\cdots,u}_{\beta})\}.$$\qed\\

According to Theorem 3.3, we can get the following
corollary immediately.\\

\noindent\textbf{Corollary 3.5.}~\emph{Let $C$ be a one-Lee weight
$\mathbb{Z}_2\mathbb{Z}_2[u]$-additive code of type $(\alpha, \beta;
k_{0}, k_{1}, k_{2})$, then its Gray image $\Phi(C)$ has parameters $[\alpha+2\beta,
k_{0}+2k_{1}+k_{2}, \lambda\frac{|C|}{2}]$}. \\

\noindent\textbf{Example 3.6.} Consider the
$\mathbb{Z}_2\mathbb{Z}_2[u]$-additive code $C$ of type (7, 1; 1, 0,
1) generated by the following generator matrix
$$ G=\begin {pmatrix}
 1 & 1 & 1 & 0 & 0 &0 & 1 & u \\
 0 & 0 & 0 & 1 & 1 &1 & 1 & u
\end{pmatrix}.$$ The Lee weight enumerator of $C$ is
$W_{L}(x,y)=x^9 + 3x^3y^6$. So $C$ is a one-Lee weight
$\mathbb{Z}_2\mathbb{Z}_2[u]$-additive code. The Gray image of $C$ is a binary linear quasi-cyclic code of
index 3 with good parameters [9, 2, 6], and the Gray image of $C^{\perp}$ is a binary $[9,7,2]$-linear code.\\

\noindent\textbf{Example 3.7.} Consider the
$\mathbb{Z}_2\mathbb{Z}_2[u]$-additive code $C$ of type (3, 3; 1, 0, 1)
generated by the following generator matrix
\[G=\left( \begin{matrix}
   1 & 0 & 1 & u & 0 & u  \\
   0 & 1 & 1 & u & u & 0  \\
\end{matrix} \right)\]

The code $C$ has the same Lee weight enumerator as the code in
Example 3.6. So, $C$ is also a one-Lee weight
$\mathbb{Z}_2\mathbb{Z}_2[u]$-additive code. The $\Phi(C)$ is an
optimal
binary $[9,2,6]$-linear code, however, $\Phi(C)$ is not quasi-cyclic.\\

\noindent\textbf{Example 3.8.} Consider the
$\mathbb{Z}_2\mathbb{Z}_2[u]$-additive code $C$ of type (4, 5; 2, 0,
1) generated by the following generator matrix
$$ G=\begin {pmatrix}
 1 & 1 & 1 & 1 & u & 0 & u & 0 & 0 \\
 0 & 0 & 1 & 1 & u & u & 0 & u & 0 \\
 0 & 0 & 0 & 0 & u & u & u & 0 & u \\
\end{pmatrix}.$$
The Lee weight enumerator of $C$ is $W_{L}(x,y)=x^{14} + 7x^6y^8$.
So, $C$ is a one-Lee weight
$\mathbb{Z}_2\mathbb{Z}_2[u]$-additive code.
$\Phi(C)$ is an optimal binary $[14,3,8]$-linear code, and $\Phi(C^{\perp})$ is a binary $[14,11,2]$-linear code.\\

\dse{4~~One-Lee weight $\mathbb{Z}_2\mathbb{Z}_2[u]$-additive formally self-dual codes}

 Let $C$ be a one-Lee weight
$\mathbb{Z}_2\mathbb{Z}_2[u]$-additive code with weight $m$. Let
$N=\alpha+2\beta$ and $\{A_{0},A_{1},\cdots,A_{N}\}$ be its Lee
weight distribution. The Lee weight enumerator of $C$ is defined by
$$W_{L}(x,y)=\sum_{i=0}^{N}A_{i}x^{N-i}y^{i}=x^{N}+(|C|-1)x^{N-m}y^{m}.$$
The MacWilliams identity with respect to the Lee weight enumerator of $C$ is
$$W_{L_{C}}(x,y)=\frac{1}{|C^{\perp}|}W_{L_{C^{\perp}}}(x+y,x-y).$$
Let $\{B_{0},B_{1},\cdots,B_{N}\}$ be the Lee weight distribution of
$C^{\perp}$. Then we have
\begin{eqnarray}
x^{N}+(|C|-1)x^{N-m}y^{m}=\frac{1}{|C^{\perp}|}\sum_{i=0}^{N}B_{i}(x+y)^{N-i}(x-y)^{i}.
\end{eqnarray}
Let $x=1$, then by differentiating (1) with respect to $y$ and
setting $y=1$, we obtain
\begin{eqnarray}
|C|=\frac{1}{|C^{\perp}|} 2^{N},
\end{eqnarray}
\begin{eqnarray}
m(|C|-1)=\frac{|C|}{2}(N-B_{1}),
\end{eqnarray}
\begin{eqnarray}
m^{2}(|C|-1)=\frac{|C|}{2}[\frac{N}{2}(N+1)-NB_{1}+B_{2}].
\end{eqnarray}

\noindent\textbf{Theorem 4.1.}\emph{ Let $C$ be a one-Lee weight
$\mathbb{Z}_2\mathbb{Z}_2[u]$-additive code with weight $m$, then
$d_{L}(C^{\perp})\geq2.$ If $\lambda=1$, then $d_{L}(C^{\perp})\geq3$. If
$\lambda=2$ and $C^{\perp}$ is a one-Lee weight additive code,
then $C=C^{\perp}$ and $C=\{(0,0|),(1,1|)\}$ or $C=\{(|0),(|u)\}.$}\\

\noindent\pf Since $m=\lambda\frac{|C|}{2}$ and
$N=\alpha+2\beta=\lambda(|C|-1)$, then
$m(|C|-1)=\lambda\frac{|C|}{2}(|C|-1)=\frac{|C|}{2}N.$ By (3), we
have $m(|C|-1)=\frac{|C|}{2}(N-B_{1})$, therefore we can deduce that
$B_{1}=0$ and $d(C^{\perp})\geq2.$ By (4), we have
$m^{2}(|C|-1)=\frac{|C|}{2}[\frac{N}{2}(N+1)-NB_{1}+B_{2}].$ So we
obtain that $B_{2}=\frac{N}{2}(\lambda-1)$. Clearly $\lambda=1$, so
$B_{2}=0$ and $d(C^{\perp})\geq3$.

If  $\lambda=2$, then
$B_{2}=\frac{N}{2}(\lambda-1)=\frac{N}{2}=|C|-1.$  Since $C^{\perp}$
is a one-Lee weight additive code, then
$|C^{\perp}|=B_{0}+B_{2}=|C|$. By (2), $|C|=\frac{1}{|C^{\perp}|}
2^{N}$. Therefore, we can deduce that $m=|C|=N=2.$ So $\alpha=2,
\beta=0$ or $\alpha=0, \beta=1$. If $\alpha=2$ and $\beta=0$,
then $(1,1|)$ is the only codeword with weight 2.
If $\alpha=0$ and $\beta=1$, then $(|u)$ is the only codeword with weight 2.
So $C=C^{\perp}$, and $C=\{(0,0|),(1,1|)\}$ or $C=\{(|0),(|u)\}.$\qed\\

A linear code $C$ is called self-orthogonal if $C\subseteq
C^{\perp}$ and self-dual if $C= C^{\perp}$. A linear code $C$ is
called formally self-dual if its weight enumerator is the same as
the weight enumerator of its orthogonal. Clearly any self-dual code
is necessarily formally self-dual,
but not all formally self-dual code are self-dual.\\

\noindent\textbf{Definition 4.2.} Let $C$ be a
$\mathbb{Z}_2\mathbb{Z}_2[u]$-additive codes.
If $W_{L_{C}}(x,y)=W_{L_{C^{\perp}}}(x,y)$, then $C$ is called a $\mathbb{Z}_2\mathbb{Z}_2[u]$-additive formally self-dual code.\\

\noindent\textbf{Example 4.3.} Consider the code
$C=\{(0,0|),(1,1|)\}$ in Theorem 4.1. It is clear that $C=C^{\perp}$
and their weight enumerators are both $W_{L_{C}}(x,y)=x^{2}+y^{2}$.
Hence, $C$ is a formally self-dual code. But if the code
$C=\{(0,0|),(1,0|)\},$ then it is clear that
$C^{\perp}=\{(0,0|),(0,1|)\}$
and their weight enumerators are both $W_{L_{C}}(x,y)=x^{2}+xy$. Hence, $C$ is a formally self-dual code, but not self-dual.\\

Now, we give the classification of the one-Lee weight
$\mathbb{Z}_2\mathbb{Z}_2[u]$-additive formally self-dual codes  and
gives possible weights for all one-Lee weight formally self-dual
codes. It is clear that if $C$ is a one-Lee weight
$\mathbb{Z}_2\mathbb{Z}_2[u]$-additive
formally self-dual code, then $\Phi(C)$ is a one weight binary formally self-dual code.\\

\noindent\textbf{Theorem 4.4.} \emph{Let $C$ be a one-Lee weight
formally self-dual code in
$\mathbb{Z}_2^{\alpha}R^{\beta}$.
Then the weight of $C$ is even if and only if $(1_{\alpha}|u_{\beta})\in C$ and $\alpha$ is even.}\\

\noindent\pf If $C$ is a one-Lee weight formally self-dual code in
$\mathbb{Z}_2^{\alpha}R^{\beta}$, then $C$ and
$C^{\perp}$ have the same  weight enumerators. Hence, $C^{\perp}$ is
also a one-Lee weight code.

For any
$c=(a|b)=({a_{0},a_{1},\cdots,a_{\alpha-1}|b_{0},b_{1},\cdots,
b_{\beta-1}})\in C$,
\begin{align*}
(a|b)\cdot(1_{\alpha}|u_{\beta})&= u\sum_{i=0}^{\alpha-1}a_{i}+u\sum_{j=0}^{\beta-1}b_{j}\\
&= u\sum_{i=0}^{\alpha-1}w_{L}(a_{i})+u\sum_{j=0}^{\beta-1}w_{L}(b_{j})= uw_{L}(a)+uw_{L}(b)\\
&=uw_{L}(c)\in \mathbb{Z}_2+u\mathbb{Z}_2.
\end{align*}
Therefore, if the weight of $C$ is even, then
$(a|b)\cdot(1_{\alpha}|u_{\beta})=0$, so we have
$(1_{\alpha}|u_{\beta})\in C^{\perp}$. Otherwise,
$(a|b)\cdot(1_{\alpha}|u_{\beta})=u$, then
$(1_{\alpha}|u_{\beta})\notin C^{\perp}$. Note that the weight
enumerators of $C$ and $C^{\perp}$ are the same and
$(1_{\alpha}|u_{\beta})$ is the only vector with Lee weight
$\alpha+2\beta$. If $(1_{\alpha}|u_{\beta})\in C^{\perp}$, then
$(1_{\alpha}|u_{\beta})\in C$. Since the weight of $C$ is even and equals $\alpha+2\beta$, it follows that $\alpha$ is even.\qed\\

Suppose that $C$ is a one-Lee weight
 formally self-dual
code in $\mathbb{Z}_2^{\alpha}R^{\beta}$ with weight
$m$. Then the weight enumerators of $C$ and $C^{\perp}$ are the same
and the weight enumerators are
$$W_{L_{C}}(x,y)=W_{L_{C^{\perp}}}(x,y)=x^{N}+(|C|-1)x^{N-m}y^{m}.$$
By the MacWilliams identity we have
$W_{L_{C}}(x,y)=\frac{1}{|C^{\perp}|}W_{L_{C^{\perp}}}(x+y,x-y).$
Then, we have
\begin{eqnarray}
x^{N}+(|C|-1)x^{N-m}y^{m}=\frac{1}{|C^{\perp}|}((x+y)^{N}+(|C^{\perp}|-1)(x+y)^{N-m}(x-y)^{m}).
\end{eqnarray}
Note that $C$ is a formally self-dual code, so $|C|=|C^{\perp}|$.
Next, by studying the equation (5), we determine conditions for the
values of $\alpha$ and $\beta$ and give a complete classification
for this class of one-Lee weight formally self-dual code. We
consider the following two possibilities
\begin{itemize}
\item[(i)] Case 1: $N=m$. By (5), we have
\begin{eqnarray}
 x^{N}+(|C|-1)y^{N}=\frac{1}{|C|}((x+y)^{N}+(|C|-1)(x-y)^{N}).
\end{eqnarray}
Comparing the coefficients of $y^{N}$ on the both side of (6), we
can get
$$|C|-1=\frac{1}{|C|}(1+(|C|-1)(-1)^{N}).$$
If $N$ is even, we can deduce $|C|=2$. Noting that
$|C|=\frac{1}{|C^{\perp}|} 2^{N}$ and $|C|=|C^{\perp}|$, we have
$N=m=2$. In this case, it follows from Theorem 4.1 that
$C=C^{\perp}=\{(0,0|),(1,1|)\}$ or $C=C^{\perp}=\{(|0),(|u)\}.$ If
$N$ is odd, we have $|C|-1=\frac{1}{|C|}(2-|C|)$. Since $|C|$ is a
positive integer, the equation has no integer root. Therefore, $N$
is not odd.
\item[(ii)] Case 2: $N\neq m$. Clearly, $N>m$. Then, comparing the coefficient of $y^{N}$ in both sides of (5), we have
\begin{eqnarray}
0=\frac{1}{|C|}(1+(|C|-1)(-1)^{m})
\end{eqnarray}
If $m$ is even, then the right side of (7) is equal to 1. This gives a contradiction.\\
If $m$ is odd, we have $0=\frac{1}{|C|}(2-|C|)$. Hence, $|C|=2.$
Note that $|C|=\frac{1}{|C^{\perp}|} 2^{N}$ and $|C|=|C^{\perp}|$.
Therefore, $N=2.$ Since $N=\alpha+2\beta$, it follows that
$\alpha=2, \beta=0$ or $\alpha=0, \beta=1$. On the other hand,
$N>m$. So, $m=1$. If $\alpha=2$ and $\beta=0$, the $(0,1|)$ or
$(1,0|)$ are the codewords with weight 1. If $\alpha=0$ and
$\beta=1$, then there exist no formally self-dual code
with weight 1. So, $C=\{(0,0|),(0,1|)\}$ and $C^{\perp}=\{(0,0|),(1,0|)\},$ or $C=\{(0,0|),(1,0|)\}$ and $C^{\perp}=\{(0,0|),(0,1|)\}.$\\
\end{itemize}

From the discussions above, we have the following theorem.\\

\noindent\textbf{Theorem 4.5.}\emph{ Let $C$ be a one-Lee weight
$\mathbb{Z}_2^{\alpha}R^{\beta}$ formally self-dual
code with weight $m$.
\begin{itemize}
\item[\emph{(1)}] If $m$ is even, then $m=2$ and $C=C^{\perp}=\{(0,0|),(1,1|)\}$ or $C=C^{\perp}=\{(|0),(|u)\}.$
\item[\emph{(2)}] If $m$ is odd, then $m=1$ and $C=\{(0,0|),(0,1|)\}$ and $C^{\perp}=\{(0,0|),(1,0|)\},$ or $C=\{(0,0|),(1,0|)\}$ and $C^{\perp}=\{(0,0|),(0,1|)\}.$
\end{itemize}}

\dse{5~~Two-Lee weight projective
$\mathbb{Z}_2\mathbb{Z}_2[u]$-additive codes}

In this section, we concentrate our study on two-Lee weight
projective $\mathbb{Z}_2\mathbb{Z}_2[u]$-additive codes and
determine the algebraic structure of two-Lee weight projective code.
Recall that a code is said to be $t$-Lee weight code if the
cardinality of the set of nonzero Lee weights is $t$.
A code over $\mathbb{Z}_2\mathbb{Z}_2[u]$ is said to be projective code if the nonzero minimum Lee weight of its dual code is at least three.\\

\noindent\textbf{Definition 5.1.} A $\mathbb{Z}_2\mathbb{Z}_2[u]$-additive code is called a
two-Lee weight projective $\mathbb{Z}_2\mathbb{Z}_2[u]$-additive code if
the cardinality of the set of nonzero Lee weight is two and the minimum Lee weight of its dual code is at least three.\\

 Let $C$ be a $t$-Lee weight $\mathbb{Z}_2\mathbb{Z}_2[u]$-additive code, then $\Phi(C)$ is a $t$-weight binary code. Let $C$ be a two-Lee weight projective
$\mathbb{Z}_2\mathbb{Z}_2[u]$-additive code with weights $m_{1}$ and
$m_{2}$. Let $N=\alpha+2\beta$ and $\{A_{0},A_{m_{1}},A_{m_{2}}\}$
be its Lee weight distribution. The Lee weight enumerator of $C$ is
defined by
$W_{L}(x,y)=x^{N}+A_{m_{1}}x^{N-m_{1}}y^{m_{1}}+A_{m_{2}}x^{N-m_{2}}y^{m_{2}}.$
The MacWilliams identity with respect to the Lee weight enumerators
of $C$ is
$W_{L_{C}}(x,y)=\frac{1}{|C^{\perp}|}W_{L_{C^{\perp}}}(x+y,x-y).$
Let $\{B_{0},B_{1},\cdots,B_{N}\}$ be the Lee weight distribution of
$C^{\perp}$. Similar to the equations (1), (2), (3) and (4), we have
\begin{eqnarray}
A_{0}+A_{m_{1}}+A_{m_{2}}=|C|=\frac{1}{|C^{\perp}|} 2^{N},
\end{eqnarray}
\begin{eqnarray}
m_{1}A_{m_{1}}+m_{2}A_{m_{2}}=\frac{|C|}{2}(N-B_{1}),
\end{eqnarray}
\begin{eqnarray}
m_{1}^{2}A_{m_{1}}+m_{2}^{2}A_{m_{2}}=\frac{|C|}{2}[\frac{N}{2}(N+1)-NB_{1}+B_{2}].
\end{eqnarray}

\noindent\textbf{Theorem 5.2.}\emph{ Let $C$ be a two-Lee weight
projective $\mathbb{Z}_2\mathbb{Z}_2[u]$-additive code with weight
$m_{1}$ and $m_{2}$. Assume $N=\alpha+2\beta$, then
\begin{eqnarray}
N^{2}-N(2m_{1}+2m_{2}-1)+m_{1}m_{2}\left(4-\frac{4}{|C|}\right)=0,
\end{eqnarray}
\begin{eqnarray}
A_{m_{1}}=\frac{\frac{|C|}{2}N-m_{2}(|C|-1)}{m_{1}-m_{2}},~~~~A_{m_{2}}=\frac{\frac{|C|}{2}N-m_{1}(|C|-1)}{m_{2}-m_{1}}.
\end{eqnarray}}
\noindent\pf Since $C$ is a two-Lee weight projective
$\mathbb{Z}_2\mathbb{Z}_2[u]$-additive codes, $A_{0}=1$ and
$B_{1}=B_{2}=0.$ Let us consider the integral coefficient polynomial
$$(x-m_{1})(x-m_{2})=a_{0}+a_{1}x+x^{2},$$ where $a_{0}=m_{1}m_{2}$
and $a_{1}=-(m_{1}+m_{2}).$ Note that $a_{0}+a_{1}m_{1}+m_{1}^{2}=0$
and $a_{0}+a_{1}m_{2}+m_{2}^{2}=0.$ Then, computing the combination
$a_{0}(8)+a_{1}(9)+(10)$ we can get the equation (11). Finally, by solving the linear system of (8) and (9), we can obtain (12).\qed\\

\noindent\textbf{Theorem 5.3.} \emph{Let $C$ be a two-Lee weight
projective $\mathbb{Z}_2\mathbb{Z}_2[u]$-additive code with
$N=\alpha+2\beta$. Then there exists a two-Lee weight projective
$\mathbb{Z}_2\mathbb{Z}_2[u]$-additive code with nonzero weights
$\frac{N}{2}$ and $\frac{|C|}{2}$.}\\

 \noindent\pf By (11),
$\Delta=(2m_{1}+2m_{2}-1)^{2}-4m_{1}m_{2}(4-\frac{4}{|C|}).$ Suppose
that $m_{1}m_{2}=\omega|C|$ for some $\omega\in \mathbb{N}$, then
$$\Delta=(2m_{1}+2m_{2}-1)^{2}-16(m_{1}m_{2}-\omega).$$ Since $N$ is
an integer, $\Delta$ must be a square. Hence, when
$\omega=\frac{m_{2}}{2}$, we have
\begin{eqnarray}
\Delta=(2m_{1}-2m_{2}-1)^{2},~~~m_{1}=\frac{|C|}{2}.
\end{eqnarray}
Thus
\begin{eqnarray}
N=\frac{(2m_{1}+2m_{2}-1)\pm (2m_{1}-2m_{2}-1) }{2}.
\end{eqnarray}
By (14), $N=2m_{1}-1$ or $N=2m_{2}$. When $N=2m_{1}-1$, then
$m_{1}=\frac{N+1}{2}$. In (13), $m_{1}=\frac{|C|}{2}$. So $N+1=|C|$.
In (12),
$A_{m_{2}}=\frac{\frac{|C|}{2}N-m_{1}(|C|-1)}{m_{2}-m_{1}}$, so we
can deduce $A_{m_{2}}=0$, this contradicts the definition of the
two-Lee weight codes. When $N=2m_{2}$, then $m_{2}=\frac{N}{2}$. \qed\\

Next, using the above properties, we construct some two-Lee weight
(projective) $\mathbb{Z}_{2}\mathbb{Z}_{2}[u]$-additive
codes which produce optimal binary linear codes.  \\

\noindent\textbf{Example 5.4.} Consider the
$\mathbb{Z}_{2}\mathbb{Z}_{2}[u]$-additive code $C$ with $\alpha=8$
and $\beta=4$ generated by the following generator matrix
\[G=\left(
  \begin{array}{cccccccccccc}
     0 & 0 & 0 & 0 & 1 & 1 & 1 & 1 & 1 & 1 & 1+u & 1+u  \\
   1 & 1 & 0 & 1 & 1 & 0 & 1 & 0 & 0 & 0 & u & u  \\
   1 & 0 & 1 & 1 & 0 & 0 & 1 & 1 & 0 & u & 0 & u  \\
   1 & 0 & 0 & 0 & 0 & 0 & 0 & 0 & u & u & u & u  \\
  \end{array}
\right)\]

The code $C$ is of type (8, 4; 2, 0, 2) and the Lee weight
enumerator is $W_{L_{C}}(x,y)=x^{16} + 7x^8y^8 + 8x^7y^9.$ So, the
code is a two-Lee weight $\mathbb{Z}_2\mathbb{Z}_2[u]$-additive
code. The Gray image of $C$ is an optimal binary $[16,4,8]$-linear code.\\

\noindent\textbf{Example 5.5.} Consider the
$\mathbb{Z}_{2}\mathbb{Z}_{2}[u]$-additive code $C$ with $\alpha=6$
and $\beta=4$ generated by the following generator matrix
$$ G=\begin{pmatrix}
 1 & 0 & 0 & 1 & 0 & 1 & 1 & 1 & 1 & 1  \\
   0 & 1 & 0 & 0 & 1 & 1 & 0 & u & 0 & u  \\
   0 & 0 & 1 & 1 & 1 & 1 & 0 & 0 & u & u  \\
   0 & 0 & 0 & 0 & 0 & 0 & u & u & u & u
\end{pmatrix}.$$

The code $C$ is of type (6, 4; 2, 0, 2) and the Lee weight
enumerator is $W_{L_{C}}(x,y)=x^{14} + 8x^7y^7 + 7x^6y^8.$ So, $C$
is a two-Lee $\mathbb{Z}_2\mathbb{Z}_2[u]$-additive code. Moreover,
$\Phi(C)$ is an optimal binary $[14,4,7]$-linear code. By the
MacWilliams identity, we have
\begin{eqnarray*}
&W_{L_{^{C^{\perp}}}}(x,y)&=x^{14} + 28x^{11}y^3 + 77x^{10}y^4 +
112x^9y^5 + 168x^8y^6 +
    232x^7y^7\\
&& + 203x^6y^8 + 112x^5y^9 + 56x^4y^{10} + 28x^3y^{11} + 7x^2y^{12}.
\end{eqnarray*}
Since the nonzero minimum Lee weight of $C^{\perp}$ is at least
three, the code $C$ is a two-Lee weight projective
$\mathbb{Z}_2\mathbb{Z}_2[u]$-additive code.
The Gray image of $C^{\perp}$ gives an optimal binary $[14,10,3]$-linear code.\\

\noindent\textbf{Example 5.6.} Consider the
$\mathbb{Z}_{2}\mathbb{Z}_{2}[u]$-additive code $C$ with $\alpha=8$
and $\beta=8$ generated by the following generator matrix
\[G=\left(
  \begin{array}{cccccccccccccccc}
     0 & 0 & 0 & 0 & 0 & 0 & 0 & 0 & u & u & u & u & u & u & u & u  \\
   1 & 1 & 1 & 1 & 1 & 1 & 1 & 1 & 0 & 0 & 0 & 0 & u & u & u & u  \\
   0 & 0 & 0 & 0 & 1 & 1 & 1 & 1 & 0 & 0 & u & u & 0 & 0 & u & u  \\
   0 & 0 & 1 & 1 & 0 & 0 & 1 & 1 & 0 & u & 0 & u & 0 & u & 0 & u  \\
   0 & 1 & 0 & 1 & 0 & 1 & 1 & 0 & 1+u & 1+u & 1+u & 1+u & 1+u & 1+u & 1+u & 1+u  \\
  \end{array}
\right)\]

The code $C$ is of type (8, 8; 3, 0, 2) and the Lee weight
enumerator is $$W_{L_{C}}(x,y)=x^{24}+28x^{12}y^{12}+3x^8y^{16}.$$
Further,
\begin{eqnarray*}
&W_{L_{^{C^{\perp}}}}(x,y)&=x^{24} + 64x^{21}y^3 + 378x^{20}y^4 +
1344x^{19}y^5+4032x^{18}y^6+ 10752x^{17}y^7\\
&& + 23439x^{16}y^8 +
    40960x^{15}y^9 + 60480x^{14}y^{10} + 77952x^{13}y^{11} + 85484x^{12}y^{12}\\
&& + 77952x^{11}y^{13} + 60480x^{10}y^{14} +
    40960x^9y^{15} + 23439x^8y^{16}+ 10752x^7y^{17}\\
&&+ 4032x^6y^{18} + 1344x^5y^{19} + 378x^4y^{20} + 64x^3y^{21} +
y^{24}.
\end{eqnarray*}
Note that the nonzero minimum Lee weight of $C^{\perp}$ is at least
three. So, the code $C$ is a two-Lee weight projective
$\mathbb{Z}_2\mathbb{Z}_2[u]$-additive code.
The Gray image $\Phi(C)$ is an optimal binary $[24,5,12]$-linear code, and $\Phi(C^{\perp})$ is an optimal binary $[24,19,3]$-linear code.\\

\noindent\textbf{Example 5.7.} Consider the
$\mathbb{Z}_{2}\mathbb{Z}_{2}[u]$-additive code $C$ with $\alpha=8$
and $\beta=4$ generated by the following generator matrix

\[G=\left(
  \begin{array}{cccccccccccc}
    1 & 1 & 1 & 1 & 1 & 1 & 1 & 1 & u & u & u & u  \\
   1 & 1 & 1 & 1 & 0 & 0 & 0 & 0 & 0 & 0 & u & u  \\
   1 & 1 & 1 & 1 & 0 & 0 & 0 & 0 & u & u & 0 & 0  \\
   0 & 0 & 1 & 1 & 0 & 0 & 1 & 1 & 0 & u & 0 & u  \\
   0 & 1 & 0 & 1 & 0 & 1 & 0 & 1 & 1 & 1 & 1 & 1  \\
  \end{array}
\right)\]

The code $C$ is of type (8, 4; 3, 0, 2) and the Lee weight
enumerator is $$W_{L_{C}}(x,y)=x^{16}+30x^{8}y^{8}+y^{16}.$$
Further,
$$W_{L_{C^{\perp}}}(x,y)=x^{16} + 140x^{12}y^4 + 448x^{10}y^6
+ 870x^8y^8 +448x^6y^{10} + 140x^4y^{12} + y^{16}.$$ Since the
nonzero minimum Lee weight of $C^{\perp}$ is at least three,  the
code $C$ is a two-Lee weight projective
$\mathbb{Z}_2\mathbb{Z}_2[u]$-additive code. The Gray image
$\Phi(C)$ is an optimal binary $[16,5,8]$-linear code, and
$\Phi(C^{\perp})$ gives an optimal binary $[16,11,4]$-linear code.

\dse{6~~Conclusion} This paper is devoted to the study of one-Lee
weight and two-Lee weight codes over
$\mathbb{Z}_{2}\mathbb{Z}_{2}[u]$, where $u^{2}=0$. The structures
and properties of these codes are obtained. We also give a complete
classification for one-Lee weight
$\mathbb{Z}_2\mathbb{Z}_2[u]$-additive formally self-dual codes.
Finally, we study the algebraic structure of two-Lee weight
projective $\mathbb{Z}_{2}\mathbb{Z}_{2}[u]$-additive codes. By
using these structures, we obtain some optimal binary linear codes
from one-Lee weight and two-Lee weight
$\mathbb{Z}_{2}\mathbb{Z}_{2}[u]$-additive codes. This shows that
one-Lee weight and two-Lee weight
$\mathbb{Z}_{2}\mathbb{Z}_{2}[u]$-additive codes are a good resource
for constructing optimal binary linear codes. A natural problem is
to study $N$-Lee weight codes over
$\mathbb{Z}_{q}\mathbb{Z}_{q}[u]$, where $u^{2}=0$ and $q$ is a
prime power. Meanwhile, it would be interesting to classify two-Lee
weight projective $\mathbb{Z}_{2}\mathbb{Z}_{2}[u]$-additive
formally self-dual codes.\\

\end{document}